\documentclass[11pt]{article}
\usepackage{amsmath,amssymb,amsfonts}

\newtheorem{theorem}{Theorem}
\newtheorem{proposition}{Proposition}
\newtheorem{corollary}{Corollary}

\def\D{{\cal D}}
\def\R{{\mathbb R}}

\def\nil{{\mathrm{Nil}\, }}

\begin{document}

\title{Surfaces of revolution in the Heisenberg group and the spectral generalization of the Willmore functional
\thanks{The work was supported by RFBR (no. 06-01-00094a) and the
complex integration project 1.1 of SB RAS. he second author
(I.A.T.) was also supported by the Program of Basic Researches of
Ministry of Education and Science of Kazakhstan (no. F03969-4).}}
\author{D.A. Berdinsky
\thanks{Institute of Mathematics,
630090 Novosibirsk, Russia; e-mail: berdinsky@ngs.ru .} \and I.A.
Taimanov
\thanks{Institute of Mathematics, 630090 Novosibirsk, Russia;
e-mail: taimanov@math.nsc.ru.} }
\date{}
\maketitle

\section{Introduction and main results}
\label{sec1}

In this paper we study the spectral generalization of the Willmore
functional for the surfaces in the three-dimensional nilpotent Lie
group $\nil$ with left-invariant metric admitting four-dimensional
isometry group, i.e. endowed by one of Thurston's geometries.

The Weierstrass representation for surfaces in $\nil$ was
introduced by us in \cite{BT} where following the spectral point
of view on the Willmore functional adopted in \cite{T1,T2} we
proposed its generalization as
$$
E(M) = \int_M UV \frac{i dz \wedge d\bar{z}}{2}
$$
where $U$ and $V$ are the potentials of the Dirac operator coming
into the representation. For the case of surfaces in $\R^3$ this
formula gives the quarter of the Willmore functional ${\cal W} =
\int H^2 d\mu$. However for surfaces in $\nil$ the functional $E$
is not proportional to the Willmore functional which in general is
equal to $\int (H^2 + \widehat{K}) d\mu$ where $\widehat{K}$ is
the sectional curvature of the ambient space along the tangent
plane to the surface.

In this paper we demonstrate that for surfaces in $\nil$ the
functional $E(M)$ resembles the Willmore functional for surfaces
in $\R^3$ in many geometrical respects.

In particular, we prove that for closed surfaces of revolution $E$
is positive and for spheres of revolution the minima of $E$ are
given by constant mean curvature spheres (see Theorem 2 and
corollaries therein). Moreover these spheres are critical points
of $E$ (see Theorem 3).

We observe the relation of the functionals $E$ and ${\cal W}$ to
the isoperimetric problem: in particular, both functionals $E$ and
$\frac{1}{4}{\cal W}$ attain the same value $\pi$ on the cmc
spheres in $\nil$ and $\R^3$ respectively (see Theorem 1). For
$\R^3$ these spheres are isoperimetric surfaces and it is
conjectured that the same is true for $\nil$. Therewith we also
show in \S \ref{sec3} how to derive some results from
\cite{abresch} and \cite{FMP} by using the Weierstrass
representation.

We find this relation of the theory of the Willmore functional to
the isoperimetric problem interesting. We discuss this relation
and some open questions in \S \ref{sec5} and demonstrate it one
more time for the case of surfaces in $S^2 \times \R$ in \S
\ref{sec6.2}.

In \S \ref{subsec6.1} we also derive the Euler--Lagrange equation
for the functional $E$.

\section{The Weierstrass representation of surfaces in the
Heisenberg group} \label{sec2}

The Heisenberg group $\nil$ is the nilpotent Lie
group formed by all matrices of the form
$$
\left(
\begin{array}{ccc}
1 & x & z \\
0 & 1 & y \\
0 & 0 & 1
\end{array}
\right), \ \ \ \ x,y,z \in \R,
$$
with the standard multiplication. It is assumed that the group is
endowed by the left invariant metric
$$
ds^2 = dx^2 + dy^2 + (dz - x dy)^2.
$$
The Lie algebra is spanned by three generators $e_1 = e_x,
e_2=e_y, e_3=e_z$ which meet the commutativity relations
$$
[e_1,e_2] = e_3, \ \ [e_1,e_3]=[e_2,e_3]=0.
$$
The scalar product at the Lie algebra to $\nil$ induced by this
left-invariant metric we denote by
$$
\langle u,v \rangle = \sum_{i=1}^3 u^iv^i, \ \ u=\sum u^ie_i, \ v=
\sum v^k e_k.
$$
As a smooth manifold this group is diffeomorphic to $\R^3$ and on
$\nil$ we also may introduce the cylindrical coordinates
$(\rho,\phi,h)$ as follows:
$$
x = \rho \cos \phi, \ \ y = \rho \sin \phi, \ \
  z = \frac{\rho^2}{2} \cos \phi \sin \phi  +  h.
$$
Given a point $z=h$ on the $z$-axis we draw the geodesic $\gamma$
of the length $\rho$ orthogonally to the $z$-axis in the direction
defined by $\phi$, the angle between $\gamma$ and the $x$-axis.
The end-point of the geodesic has the coordinates $(\rho,\phi,h)$.

The metric in the cylindrical coordinates takes the form
\begin{equation}
\label{cylmetric} ds^2 = d \rho^2  -  \rho^2 dh d\phi  +
\frac{1}{4} \rho^2 (4 + \rho^2) d \phi^2  +  dh^2
\end{equation}
and we see that it is invariant under rotations $\phi \to \phi +
\theta$ around the $z$-axis.

In fact from this formula for the metric it is easily derived that
$\nil$ has a four-dimensional isometry group generated by left
translations: $g \to hg, h \in \nil$, and rotations around the
$z$-axis.

Let us expose basic facts on the Weierstrass representation of
surfaces in three-dimensional Lie groups $G$ introduced in
\cite{BT}.

The Weierstrass representation of a surface
$$
f: M \to G
$$
defines it in terms of a solution to a nonlinear equation
\begin{equation}
\label{nil-weingarten}
{\cal D}_\nil \psi = \left[\left(
\begin{array}{cc}
0 & \partial \\
-\bar{\partial} & 0
\end{array}
\right)
 +
\left(
\begin{array}{cc}
U & 0 \\
0 & V
\end{array}
\right) \right] \psi = 0,
\end{equation}
where $z$ is a conformal parameter on the surface,
$$
Z_1 = \frac{i}{2} ( \bar{\psi}_2^2 + \psi_1^2), \ \ \ Z_2 =
\frac{1}{2} ( \bar{\psi}_2^2 - \psi_1^2), \ \ \ Z_3 = \psi_1
\bar{\psi}_2
$$
and
$$
f^{-1} f_z = \sum_{k=1}^3 Z_k e_k
$$
is the linear expansion of $f^{-1}f: M \to T_1\nil$ in the
generators $e_1,e_2,e_3$ of the Lie algebra, of $\nil$, identified
with the tangent space to $\nil$ at the unity.

The nonlinearity is hidden in the potentials $U$ and $V$ and for
$G=\nil$ we have
$$
U_\nil = V_\nil =
\frac{H}{2}(|\psi_1|^2+|\psi_2|^2) +
\frac{i}{4}(|\psi_2|^2-|\psi_1|^2),
$$
where $H$ is the mean curvature.

Therewith the induced metric equals
$$
ds^2 = \left(|\psi_1|^2+|\psi_2|^2\right)^2 dz d\bar{z},
$$
the Hopf differential $A = \langle \nabla_{f_z} f_z, n \rangle
(dz)^2$ takes the form
$$
A =  (\bar{\psi}_2 \partial \psi_1 -
\psi_1
\partial \bar{\psi}_2) + i \psi_1^2 \bar{\psi}_2^2.
$$
Let us denote by $n$ the normal vector, to the surface, translated
to $T_1\nil$ by the left multiplication by $f^{-1}$. It is equal
to
\begin{equation}
\label{nil-normal}
n = e^{-\alpha}\left[ i(\psi_1
   \psi_2 - \bar{\psi}_1\bar{\psi}_2) e_1 - (\psi_1 \psi_2 + \bar{\psi}_1\bar{\psi}_2) e_2 +
(|\psi_2|^2-|\psi_1|^2) e_3)\right],
\end{equation}

The derivational equations express derivatives of $\psi$ in $z$
and $\bar{z}$ and this system is formed by (\ref{nil-weingarten})
and the equations
$$
\partial \psi_1 = \alpha_z \psi_1 + Ae^{-\alpha} \psi_2 -
\frac{i}{2} \psi_1^2 \bar{\psi}_2,
$$
$$
\bar{\partial} \psi_2 = -\bar{A} e^{-\alpha} \psi_1  +
\alpha_{\bar{z}}\psi_2-
 \frac{i}{2} \bar{\psi}_1 \psi_2^2.
$$

The derivational equations are obtained by simple straightforward
computations and one of the immediate consequences is as follows

\begin{itemize}
\item
{\sl a surface (in $\nil$) has constant mean curvature if and only
if the quadratic differential
\begin{equation}
\label{hopfdif}
\widetilde{A} dz^2 = \left (A + \frac{{Z_3}^2}{2H+i} \right)dz^2
\end{equation}
is holomorphic}. \footnote{ Since it was proved earlier in
\cite{AR} that for constant mean curvature (cmc) surfaces in the
products $S^2 \times \R$ and $H^2 \times \R$ some generalizations
of the Hopf differential are holomorphic and the same was
announced for surfaces $\nil$ and other three-spaces with
four-dimensional isometry group in \cite{abresch} and since these
results by Abresch and Rosenberg motivated us to prove the same by
our methods, in \cite{BT} we attributed the statement on the
holomorphicity of this differential for cmc surfaces to Abresch.
However the detailed analysis of the formulas from
\cite{abresch,AR} shows that the Abresch--Rosenberg differential
has the form
$$
(H+i\tau) \widetilde{A} dz^2.
$$
Therewith $\nil$ is locally considered as a one-dimensional
fibration over the flat two-space with the bundle curvature
$\tau$.}
\end{itemize}

For surfaces in $\R^3$ it is clear that Hopf differential is
holomorphic if and only if the surface has constant mean
curvature. We did manage to generalize this fact for surfaces in
$\nil$ and did not succeed for surfaces in $\widetilde{SL_2(\R)}$.
It appears to be impossible. Recently it was showed by Fern\'andez
and Mira that for surfaces in this group there are non-compact not
cmc surfaces which the differential $\widetilde{A} dz^2$ is
holomorphic however all compact surfaces for which this
differential is holomorphic are cmc surfaces (see this and similar
results related to surfaces in other three-spaces with
four-dimensional isometry group in \cite{FM}).

In \cite{BT} we also introduce the (spinor) energy of a compact
surface $M$ without boundary in $\nil$ (and also in
$\widetilde{SL_2(\R)})$ as
$$
E(M) = \int_M UV \frac{i dz \wedge d\bar{z}}{2}.
$$
For a surface in $\R^3$ $U=\bar{U}=V$ in its Weierstrass
representation and $E(M) = \frac{1}{4}{\cal W}$ where ${\cal W}$
is the Willmore functional \cite{T1}. The point of view based on
the spectral theory of Dirac operators $\D$ coming in the
Weierstrass representations and taken and demonstrated in
\cite{T1,T2} assumes that the spectral properties of $\D$ have to
have important geometrical meanings. Thus we treat the functional
$E$ as the spectral generalization of the Willmore functional.
Although the product $UV$ is complex-valued it was showed in
\cite{BT} that the integral taken over a compact surface without
boundary equals
\begin{itemize}
\item
\begin{equation}
\label{energy} E(M) = \frac{1}{4} \int_M \left( H^2 +
\frac{\widehat{K}}{4} - \frac{1}{16} \right) d\mu
\end{equation}
{\sl where $\widehat{K}$ is the sectional curvature of $\nil$
along the tangent plane to a surface and $d\mu = e^{2\alpha}dx
\wedge dy$ is the induced measure on $M$.}
\end{itemize}

\section{Constant mean curvature spheres in the Heisenberg group}
\label{sec3}

\subsection{The main equation for surfaces with $\widetilde{A}=0$}

Let us first formulate some simple identities obtained from the
derivational formulas and checked by straightforward
computations:
\begin{equation}
\label{rem} \frac{\partial n_3}{\partial z} = \left(-H +
\frac{i}{2}\right) Z_3  - 2 e^{ -2\alpha} A \bar{Z}_3
\end{equation}
where $ n_3 = \langle n, e_3 \rangle$,
\begin{equation}
\label{metrics} e^{2\alpha} = \frac {4 |Z_3|^2} {1 - {n_3}^2}
\end{equation}
and
\begin{equation}
\label{eq2} \frac{\partial \overline{Z}_3}{\partial z} = ( 2H - i
) |Z_3|^2 \frac {n_3}{1-{n_3}^2}.
\end{equation}
The formula (\ref{rem}) easily follows from the derivational
equations written in terms of the immersion $f$:
$$
\nabla_{f_z} n   =  -Hf_z  -  2 A e^{-2\alpha} f_{\overline z}.
$$

Let us suppose that the differential $\widetilde{A}dz^2$ vanishes
which, in particular, implies that the mean curvature of a surface
is constant \cite{BT}:
\begin{equation}
\label{Hopf for cmc spheres} A = - \frac{{Z_3}^2}{2H+i}, \ \ \ H =
\mathrm{const}.
\end{equation}

Substituting (\ref{Hopf for cmc spheres}) into (\ref{rem}) and
expressing $e^{2\alpha}$ via (\ref{metrics}), we obtain
\begin{equation}
\label{eq1} \frac{\partial n_3}{\partial z} = \left(-H +
\frac{i}{2} + \frac {1-{n_3}^2} {4H +2i}\right) Z_3
\end{equation}
which together with (\ref{eq2}) implies
\begin{equation}
\label{main equation} \Delta  n_3 + \frac{2 n_3 ({n_{3x}}^2 +
{n_{3y}}^2)}{1 - {n_3}^2} = 0.
\end{equation}

We have

\begin{proposition}
For a surface in $\nil$ with vanishing differential
$\widetilde{A}dz^2$ the equation (\ref{main equation}) holds.

Moreover the metric $e^{2\alpha}$ on this surface is uniquely
reconstructed from the function $n_3$ and the constant $H$.
\end{proposition}

The first statement is already derived. To prove the second part
it is enough to reconstruct $Z_3$ from \eqref{eq1} and then by
using (\ref{metrics}) derive
\begin{equation}
\label{metrics2} e^{2\alpha} = \frac{4}{1-n_3^2} \frac{16H^2 +
4}{(4H^2 + {n_3}^2)^2}  \left| \frac{\partial n_3}{\partial z}
\right| ^2.
\end{equation}

\subsection{CMC spheres of revolution}

Constant mean curvature of revolution and more general cmc
surfaces with helicoidal symmetry were described in different
terms in \cite{FMP}. We demonstrate here how a description of cmc
spheres is straightforwardly derived via the Weierstrass
representation. Moreover these computations will be necessary for
us for computing the values of different functionals (area,
bounded volume, spinor energy) on these spheres.

Let us consider the following solutions to \eqref{main equation}:
\begin{equation}
\label{n3}
n_3 = \frac{r^2 - 1}{r^2 + 1},
\end{equation}
where $r^2 = x^2 + y^2, z = x+iy$.

If such a solution corresponds to a surface with $\widetilde{A}=0$
then, by \eqref{metrics2}, the induced metric
$e^{2\alpha}dzd\bar{z}$ on the surface takes the form
\begin{equation}
\label{metric3} e^{2\alpha} =
\frac{16(1 + 4H^2)(1+ r^2)^2}{(( r^2 - 1)^2 + 4 H^2(1 + r^2)^2)^2}
\end{equation}

We look for a surface which is obtained by a revolution of curve
$\gamma(r) = (\rho (r), \psi (r), h(r))$ and for which the induced
metric takes the form \eqref{metric3} where $x = r \cos \theta$ è
$y = r \sin \theta$ are the conformal coordinates on the surface with  $\theta
$ the angle of rotation.

The induced metric on a surface of revolution in  the coordinates
$r$ and $\theta$ is equal to
$$
\left( \rho^2 + \frac{\rho^4}{4} \right)d\theta^2 + \left( \left(2\rho^2 +
\frac{\rho^4}{2}\right)\psi^\prime - \rho^2 h^\prime \right)dr d\theta +
$$
$$
+ \left(h^{\prime 2} + \rho^2 \psi^{\prime 2} + \frac{1}{4} \rho^4
\psi^{\prime 2} - \rho^2 h^\prime \psi^\prime + \rho^{\prime 2} \right)
dr^2.
$$
Such a metric takes the form $e^{2\alpha}dzd\bar{z}$ if and only if
$\rho$, $h$ and $\psi$ satisfy the following equations:
$$
\rho^2 + \frac{\rho^4}{4} = r^2 e^{2 \alpha},
\ \ \ \
\left( 2 \rho^2 + \frac{\rho^4}{2}\right) \psi^\prime - \rho^2 h^\prime = 0,
$$
$$
h^{\prime 2} + \rho^{\prime 2} +\rho^2 \psi^{\prime 2}  + \frac{1}{4} \rho^4 \psi^{\prime 2} -
\rho ^2 h^\prime \psi^\prime = e^{2 \alpha}.
$$
This system is rewritten as follows:
$$
\rho = \sqrt {\sigma},
\ \ \ \ \
h^\prime = \sqrt { 1 + \frac{\sigma}{4}} \sqrt {e^{2 \alpha}  -
\frac {\sigma^{\prime 2}}{4 \sigma}},
\ \ \ \ \
\psi^\prime = \frac {1}{2} \sqrt {\frac{e^{2 \alpha} - \frac {\sigma'^2}{4 \sigma} }{1 +
\frac {\sigma}{4}}}
$$
where
$\sigma = 2 \sqrt {1 + r^2 e^{2\alpha}} -2$.

If $e^{2\alpha}$ takes the form \eqref{metric3} then
$$
\sigma = \frac{16 r^2}{ (r^2 -1)^2 + 4H^2(r^2 + 1)^2}
$$
and we have
$$
\rho = \sqrt {\frac{16 r^2}{ (r^2 -1)^2 + 4H^2(r^2 + 1)^2}},
$$
$$
h^\prime = \frac {16H(1+4H^2)r(1+r^2)^2}{((r^2-1)^2 + 4H^2(1+r^2)^2)^2},
$$
$$
\psi^\prime = 8 \frac{Hr}{(r^2-1)^2 + 4H^2(1 + r^2)^2 }.
$$
The final formulas for the generating curve $\gamma(r)$ of the surface of revolution
are as follows:
\begin{equation}
\label{cmc}
\begin{split}
\rho = \frac{ 4r }{\sqrt{(r^2-1)^2 + 4H^2 (r^2+1)^2} },\\
h = \frac{1+4H^2}{4H^2}
\left(-\frac{4H(1-r^2 + 4H^2(1+r^2))}{(r^2-1)^2 +
16H^4(1+r^2)^2+8H^2(1+r^4)} \right.+
\\
+ \left.\arctan\left[\frac{1}{4H} (r^2 -1 + 4H^2(r^2 +1))\right]\right),
\\
\psi = \arctan \left[\frac{1}{4H} (4H^2 - 1 + (1+4H^2)r^2)\right],
\end{split}
\end{equation}
where $r \in [ 0 , \infty ]$.

The following proposition is checked by straightforward computations.

\begin{proposition}
For any $H, 0 < H < \infty$, the curve \eqref{cmc}
generates by revolution a sphere with constant mean curvature $H$.
\end{proposition}

Let $T_1 \nil$ be the $S^1$-fiber bundle over $\nil$ formed by all
unit vectors. We denote by $\widehat{f}: M \rightarrow T_1 \nil$
the Gauss map which corresponds to a point $p \in M$ the unit
normal vector at $p$.

\begin{proposition}
\label{prop3}
Given $H, 0 < H < \infty$, for any point $q \in T_1
\nil$ there exists a sphere of revolution $M$  with constant mean
curvature $H$ such that $q \in \widehat{f}(M)$.
\end{proposition}

Of course, here and in the sequel we mean by spheres of revolution
not only spheres given by \eqref{cmc} but also their left
translates in the group $\nil$.

{\it Proof} of Proposition \ref{prop3}. Let $q=(p,\xi)$ with $p
\in \nil$ and $\xi \in T_p\nil$. Given the sphere $S_H$, it
follows from \eqref{cmc} that $n_3$ takes all values from $-1$
till $1$. Therefore let us take $p_1 \in S_H$ such that $\xi_3 =
n_3(p_1)$ and translate $S_H$ into a cmc sphere $S_1$ by
left-translation $g \to hg$ such that $hp_1=p$. The normal to
$S_1$ at $p_1$ equals $(\xi^\prime_1,\xi^\prime_2,\xi_3)$. Then we
rotate $S_1$ around the $z$-axis coming through $p$ to achieve a
sphere $S_2$ for which the normal at $p$ is equal to $\xi$. This
proves the proposition.

\subsection{CMC spheres}

To finish the description of all cmc spheres in $\nil$ we are left
to show all such spheres are just spheres of revolution: this fact
was proved in \cite{AR} for surfaces in $S^2 \times \R$ and $H^2
\times \R$ and was stated in \cite{abresch} for other
three-manifolds with four-dimensional isometry group. Moreover in
\cite{abresch} it is explained that the proof for the latter case
is almost the same as for the cases of products \cite{AR}. Here we
expose such a proof for in the particular case of $\nil$.

\begin{proposition}
[Abresch--Rosenberg \cite{abresch}] Given $H, 0 < H < \infty$, any
complete surface with $\widetilde{A}=0$ is a sphere of revolution.
\end{proposition}

{\it Proof.} One of the Gauss--Weingarten equations reads
\begin{equation}
\label{conformal} \nabla _ {f_z} n = - H f_{z} -
2Ae^{-2\alpha}f_{\overline{z}}.
\end{equation}
Since $\widetilde{A}=0$, we have  $A = - \frac {Z_3 ^2}{2H + i}$
and therefore at any point $p$ the vectors $\nabla_{f_x} n$ are
$\nabla_{f_y} n$ uniquely defined by $f_z$, $f_{\overline{z}}$,
and the point $p$. Moreover we have
$$
(\nabla_{f_x} n)^i = \frac{\partial n^i}{\partial x} + \Gamma ^
{i} _{jk} (p) f^j_x n^k.
$$
The equation  \eqref{conformal} takes the same form for any
conformal coordinate $w=w(z)$ on the surface. In fact it defines a
two-plane $\Pi_q$ in $T_q\left(T_1\nil\right)$ with $q=(p,n)$ such
that $\Pi_q$ is tangent to the image of the Gauss mapping of any
surface with $\widetilde{A}=0$. Thus we have a two-dimensional
distribution $\Pi$ on $T_1\nil$. Any integral surface of this
distribution is uniquely determined by any its point and, since
through any point of $T_1 \nil$ goes the image of the Gauss map of
a cmc sphere, we conclude that all complete surfaces with
$\widetilde{A}=0$ are cmc spheres. Proposition is proved.

\subsection{Remark on the isoperimetric problem for $\nil$}

Since the metric on the sphere equals
$$
e^{2\alpha}(dr^2 + r^2d\theta^2),
$$
the area element is equal
$$
d \mu = r e^{2\alpha} dr d\theta.
$$

Substituting \eqref{metric3} into this formula we compute the area
$A(H)$ of the sphere $S_H$ with constant mean curvature $H$:
$$
A(H) = \int_{S_{H}} d\mu = 2 \pi \int_{0}^{\infty} r e^{2\alpha} dr =
$$
$$
2 \pi \int_{0}^{\infty} \frac{16(1+4H^2)r(1+r^2)^2}{((r^2-1)^2
+ 4H^2(r^2+1)^2)^2} \, dr =
$$
$$
2\pi \left( \frac{1}{H^2}+\frac{1+4H^2}{4H^3}
\left (\frac{\pi}{2}-\arctan\left[\frac{4H^2-1}{4H}\right]\right) \right).
$$

In the domain $D_H$ bounded by the sphere $S_H$ we take for coordinates the parameters
$\delta \in [0,1] $, $r \in (0,\infty) $ and $\theta
\in [0, 2 \pi]$ such that the cylindrical coordinates of the point $(\delta,r,\theta)$
are equal to
$\rho = \delta \rho(r)$, $\phi = \psi(r) + \theta$ and $h=h(r)$
where the functions $\rho(r)$, $h(r)$, and $\psi(r)$ from
\eqref{cmc} define a sphere of revolution.

We have
$$
d\rho = \delta \rho^\prime(r) dr + \rho(r) d \delta,\ \
d\phi = \psi^\prime(r) dr + d \theta, \ \
dh = h^\prime(r) dr
 $$
and, substituting these formulas to \eqref{cylmetric}, we compute the induced metric:
$$
ds^2  =  \rho^2 d\delta^2 +
2\delta\rho\rho^\prime dr d\delta +
$$
$$
(\delta^2\rho^{\prime 2} -
\delta^2\rho^2 h^\prime \psi^\prime + \delta^2 \rho^2 \psi^{\prime 2}  + \frac{1}{4} \delta^4 \rho^4
 \psi^{\prime 2} + h'^2)dr^2+
$$
$$
+ \frac{1}{4} \delta^2\rho^2 (4+\delta^2
\rho^2)d\theta^2 + (-\delta^2 \rho^2 h^\prime + 2\delta^2 \rho^2 \psi^\prime + \frac{1}{2}
\delta^4 \rho^4 \psi^\prime)dr\,d\theta,
$$
and the volume form $d\nu$:
$$
d\nu =  (\rho^2 h^\prime \delta) \, d\delta dr d \theta =
\frac{256H (1+4H^2)r^3(1+r^2)^2 \delta}{((r^2-1)^2 +
4H^2(1+r^2)^2)^3}\, d\delta dr d \theta.
$$
Therefore the volume $V(H)$ of $D_H$ equals
$$
V(H) = \int_{D_H} d\nu = \pi \int_{0}^{\infty} \frac{256H
(1+4H^2)r^3(1+r^2)^2 }{((r^2-1)^2 + 4H^2(1+r^2)^2)^3} \, dr =
$$
$$
= \frac{\pi}{16H^4} \left(4H(4H^2+3) -
(4H^2+1)(4H^2-3)\left(\frac{\pi}{2} -\arctan\left[\frac{4H^2
-1}{4H}\right]\right)\right).
$$
Finally we obtain the relation between the area $V(H)$ of a cmc
sphere and the volume $V(H)$ of the domain bounded by this sphere:
\begin{equation}
\label{isoperimetric} V(H) = \frac{4\pi}{H}  - \frac{4H^2-3}{8H}
A(H).
\end{equation}

Conjecturally the relation between $A(H)$ and $S(H)$ gives a
solution to the isoperimetric problem for $\nil$.

For a general $n$-dimensional Riemannian manifold this problem
consists in finding a hypersurface $S$ which minimizes the
$(n-1)$-volume $V_{n-1}$ among all surfaces bounding domains of
$n$-volume $d$. This surface has constant mean curvature and is
called isoperimetric and we denote its volume by $V_{n-1}(d)$.
From the geometric measure theory it is known that for $n \leq 7$
an isoperimetric hypersurface is smooth \cite{Morgan}.

For $\R^3$ the isoperimetric surfaces are the round spheres and
$V_2(d) = (24\pi d^2)^{1/3}$. This was originally proved by
Schmidt in 1930s by the symmetrization method \cite{Schmidt}
however now it also can be derived from the Alexandrov theorem
that all embedded compact cmc surfaces without boundary in $\R^3$
are homeomorphic to spheres and the Hopf theorem that all cmc
spheres in $\R^3$ are the round spheres.

The analog of the Alexandrov theorem is not known for $\nil$.
However it is very unlikely that isoperimetric surfaces in $\nil$
are non-spherical and it is a reasonable and known conjecture that
{\sl isoperimetric surfaces in $\nil$ are homeomorphic to
spheres}. If it is true the cmc spheres $S_H, 0 < H < \infty$,
give isoperimetric surfaces for all $d, 0 < d < \infty$.
\footnote{After the posting of the first version of this paper in
the internet F. Morgan pointed out to us the paper \cite{Tomter}
where the constant mean curvature spheres of revolution  in $\nil$
are described and it is proved that for small volumes, i.e. for $H
\gg 0$, they are solutions to the isoperimetric problem.} We
remark that for a compact Riemannian manifold and for small
volumes the isoperimetric hypersurfaces are homeomorphic to a
sphere \cite{MorganJohnson}.

\section{The spectral generalization of the Willmore fun\-ctional}
\label{sec4}

For closed oriented surfaces in $\nil$ the spinor energy
functional introduced in \cite{BT} is equal to
\begin{equation}
\label{nil-energy}
\begin{split}
E(M) = \int_M UV \frac{i dz \wedge d\bar{z}}{2} =
\\
= \frac{1}{4} \int_M \left(H^2 - \frac{n_3^2}{4}\right) d\mu =
\frac{1}{4} \int_M \left( H^2 + \frac{\widehat{K}}{4} -
\frac{1}{16} \right) d\mu.
\end{split}
\end{equation}

For cmc spheres $S_H$ the spinor energy takes the form
\begin{equation}
\label{Willmore for spheres} E(S_H) = \frac{\pi}{2}
\int^{\infty}_{0} \left(H^2 - \frac{1}{4} n_3^2\right) e^{2\alpha}
r dr
\end{equation}
where $n_3$ and $e^{2\alpha}$ are given by \eqref{n3} and
\eqref{metric3}. Substituting these formulas for $n_3$ and the
metric into \eqref{Willmore for spheres}, by straightforward
computations we prove

\begin{theorem}
For all cmc spheres in $\nil$ the spinor energy is equal to $\pi$:
\begin{equation}
\label{sphere-energy} E(S_H) = \pi.
\end{equation}
\end{theorem}

Let us compute the classical Willmore functional
\begin{equation}
\label{willmore}
{\cal W}(M) = \int_M (H^2 + \widehat {K}) d \mu.
\end{equation}
for these spheres. Since for surfaces in $\nil$ we have
$\widehat{K} = \frac{1}{4} - e^{-2\alpha}(|\psi_2|^2 - |\psi_1|^2)
= \frac{1}{4} - n_3^2$, it follows from \eqref{nil-energy} and
\eqref{willmore} that $ \int_{S_H}  \widehat{K} d \mu = 16\pi -
\left(4H^2-\frac{1}{4}\right) A(H)$ and finally we derive that
$$
{\cal W}(S_H) =
10\pi  + \frac{\pi}{2H^2} - \pi \frac {(1+4H^2)(3H^2-\frac{1}{4})} 2H^3 \left(\frac{\pi}{2}
 -\arctan\left [\frac{4H^2 -1}{4H}\right] \right).
$$
Therefore we see that the Willmore functional does not take a
constant value on constant mean curvature spheres.

Let us compute the spinor energy for closed surfaces of revolution.

We have the $SO(2)$-action on $\nil$ by rotations around the
$z$-axis and the quotient space $\nil/SO(2)$
is the half-plane $u \geq 0$ with the local coordinates $u=\rho$ and $v=z$
where $\rho,\phi$, and $z$ are the cylindrical coordinates.
By \eqref{cylmetric}, there is a submersion
$$
\nil \to B = \nil/SO(2)
$$
where $\nil/SO(2)$ is endowed with the metric
$$
du^2 + \frac{4 u^2}{ 4u^2 + u^4} dv^2.
$$
Let $\gamma(s) = (u(s), v(s))$ be a smooth curve in $B$ which generates by revolution a smooth surface in $\nil$.
Here we denote by $s$ the natural parameter on $\gamma$.
Let $\sigma$ be the angle between $\gamma$ and the direction $\frac{\partial}{\partial u}$.
We have the following formulas for the tangent vector $t$ and the normal vector $n$:
$$
t = (\cos \sigma, (2u) ^{-1} \sqrt{4u^2 + u^4} \sin \sigma ), \ \ \
n = (- \sin \sigma, (2u) ^{-1} \sqrt{4u^2 + u^4} \cos \sigma ).
$$
Moreover $u,v$, and $\sigma$ satisfy the following ordinary differential equations:
\begin{equation}
\label{pedrosa}
\begin{cases}
\dot{u} = \cos \sigma \\
\dot {v} = (2u) ^{-1} \sqrt{4u^2 + u^4} \sin \sigma \\
\dot {\sigma} = 2H - u^{-1} \sin \sigma,
\end{cases}
\end{equation}
where the dot denotes the derivation in $s$. It follows from \eqref{pedrosa} that
$$
H = \frac 1 2 (\dot{\sigma} + u^{-1} \sin \sigma).
$$
These formulas for $t, n$, and $H$ were derived in \cite{FMP}.

It is easy to compute
$$
n_3 = \left\langle n , \frac{\partial}{\partial
 z} \right\rangle  =  \frac{2u}{\sqrt{4u^2 + u^4}} \cos \sigma, \ \ \ d \mu  = \frac 1 2 \sqrt {4u^2 + u^4} d\theta ds.
$$
Rewriting the functional $E$ in terms of $u$ and $\sigma$ we
compute
$$
E(M) = \frac{\pi}{4} \int_{\gamma} \left[\frac{(\dot{\sigma} +
\frac{sin \sigma}{u})^2}{4} - \frac{u^2}{4u^2 + u^4}
  cos^2 \sigma \right] \sqrt{4u^2 + u^4} \, ds =
$$
$$
= \frac{\pi}{4} \int_{\gamma} \left[\frac{(\dot{\sigma} -
\frac{\sin \sigma}{ u})^2}{4} \sqrt{4u^2 + u^4} +
\frac{\dot{\sigma} \sin \sigma}{u} \sqrt{4u^2 + u^4} - \frac{u^2
\cos^2 \sigma }{\sqrt{4u^2 + u^4}} \right] ds.
$$
Now we are left to notice that
$$
\frac{\pi}{4} \int_{\gamma} \left( \frac{\dot{\sigma} \sin
\sigma}{u} \sqrt{4u^2 + u^4} - \frac{u^2}{\sqrt{4u^2 + u^4}} \cos
^2 \sigma \right) ds =
$$
$$
= -  \frac{\pi}{4} \int_{\gamma} \left( \ddot{u}\sqrt{4+u^2} +
\frac{u}{\sqrt{4+u^2}} {\dot{u}}^2 \right) ds= -\frac{\pi}{4}
\int_{\gamma} \frac{\partial}{\partial s} (\dot{u} \sqrt {4+u^2})
ds.
$$
Thus we prove the following theorem.

\begin{theorem}
Given a closed surface $M$ in $\nil$ obtained by revolving a curve
$\gamma \subset B$ around the $z$-axis, the spinor energy of $M$
equals
\begin{equation}
\begin{split}
\label{willmore1} E(M) = \frac 1 4  \int_{\gamma} \left( H^2 -
\frac 1 4 {n_3}^2 \right)
d\mu = \\
\frac{\pi}{8} \int_{\gamma} \left(\dot {\sigma} -
\frac{\sin \sigma}{ u}\right)^2 \sqrt{4u^2 + u^4} ds - \frac{\pi}{4}
\int_{\gamma} \frac{\partial{[\dot{u}
\sqrt{4+u^2}]}}{\partial s} ds = \\
\\
\frac{\pi}{8} \int_{\gamma} \left(\dot {\sigma} -
\frac{\sin \sigma}{ u}\right)^2 \sqrt{4u^2 + u^4} ds + \frac{\pi \chi(M)}{2}
\end{split}
\end{equation}
where $\chi(M)$ is the Euler characteristic of $M$.

Moreover if $\dot{\sigma} = \frac{\sin \sigma}{u}$ everywhere on
the surface then it is a cmc sphere.
\end{theorem}

\begin{corollary}
For spheres of revolution $E(M) \geq \pi$ and the equality is
attained exactly at cmc spheres.
\end{corollary}

\begin{corollary}
For tori of revolution $E(M) > 0$.
\end{corollary}

We also have

\begin{theorem}
The cmc spheres in $\nil$ are the critical points of the spinor
energy  functional $E$.
\end{theorem}

{\it Proof}. The Euler--Lagrange equation for $E$ takes the form
$$
\Delta H + 2H(H^2-K)+2 e^{-4 \alpha} (A\bar{Z}_3 ^2 + \bar{A}
 Z_3^2) = 0
$$
(see Theorem 4 in \S \ref{subsec6.1}). The cmc spheres meets the
equations $H = \mathrm{const}$ and $A = - \frac{Z_3 ^2}{2H+i}$
which imply
$$
\Delta H = 0, \ \ \ 2H(H^2 - K) = 8 e^{-4\alpha} H |A|^2 = 8
e^{-4\alpha} H \frac{|Z_3|^4}{4H^2 +1},
$$
$$
2e^{-4\alpha} (A\bar{Z}_3 ^2 + \bar{A} Z_3^2) = - 8 e^{-4 \alpha}
H \frac{|Z_3|^4}{4H^2 +1}.
$$
It follows from these formulas that the cmc spheres in $\nil$ meet
the Euler--Lagrange equation for $E$. Theorem is proved.

\section{Final remarks and open questions}
\label{sec5}

We see that the energy functional in many geometrical respects
behaves similarly to the functional
$$
\frac{{\cal W}(M)}{4} = \frac{1}{4}\int_M H^2 d\mu
$$
for closed oriented surfaces in $\R^3$.

Indeed

\begin{enumerate}
\item
as in \eqref{sphere-energy} we have
$$
\frac{{\cal W}}{4} = \pi
$$
for all isoperimetric surfaces, i.e. the round spheres, in $\R^3$;

\item
we have
$$
\frac{{\cal W}(M)}{4} = \frac{1}{4} \int_M
\left(\left(\frac{\varkappa_1-\varkappa_2}{2}\right)^2 +
\varkappa_1\varkappa_2\right)d\mu =
$$
$$
= \frac{1}{4} \int_M
\left(\frac{\varkappa_1-\varkappa_2}{2}\right)^2 d\mu +
\frac{\pi}{2}\chi(M)
$$
where $\varkappa_1$ and $\varkappa_2$ are the principal
curvatures. The latter formula is similar to \eqref{willmore1}
however the quantities $\dot{\sigma}$ and $\frac{\sin \sigma}{u}$
are not the principal curvatures of a surface of revolution;

\item
the condition $A=0$ distinguishes among complete compact surfaces
in $\R^3$ exactly the cmc spheres which are the minima of the
Willmore functional ${\cal W}$ for surfaces. Among closed surfaces
of revolution in $\nil$ the similar condition $\widetilde{A}=0$
distinguishes exactly the cmc spheres which are the minima of the
functional $E$ (among spheres of revolution).
\end{enumerate}

These geometrical observations confirm that the functional $E$
coming from the spectral theory of the Weierstrass representation
sounds to be the right generalization of the Willmore functional
to surfaces in $\nil$.

We also have to  consider the equation $\widetilde{A}=0$ as
distinguishing generalized umbilic surfaces:  both in $\R^3$ and
in $\nil$ the complete compact ``umbilic'' surfaces are the cmc
spheres. We remark that for cmc spheres in $\nil$ only the poles,
i.e points invariant under rotation symmetry, are umbilics in the
classical sense.

From the point of view of this generalization the following open
problems are interesting for study:

\begin{enumerate}
\item
to prove that $E$ is bounded from below for each topological type
of closed oriented surfaces or even to prove that $E$ is positive;

\item
to prove that the cmc spheres are the global minima of $E$ among
spheres;

\item
to generalize the formula \eqref{willmore1} for general surfaces;

\item
to find the minima of $E$ among surfaces of fixed topological type
and, in particular, to find the substitution of the Willmore
conjecture.
\end{enumerate}

Of course it sounds interesting to consider the same questions for
surfaces in $\widetilde{SL_2(\R)}$ for which case the Weierstrass
representation and the energy functional were also derived in
\cite{BT}.

\section{Appendices}

\subsection{The Euler--Lagrange equations for the
functional $E$} \label{subsec6.1}

\begin{theorem}
\label{theorem4} Let $f: M \to \nil$ be a regular surface and $r:
M \times [0,1] \to \nil$ be its smooth variation: $r_0=f$. We
assume that $r$ is constant on the boundary of $M$ if it exists:
$r(p,t)=f(p)$. Let $ \frac{\partial r(p,t)}{\partial
t}\left|_{t=0}\right. = \varphi n$ where $n$ is the unit normal
field to $M$. Then the variation of $E$ at $t=0$ equals
\begin{equation}
\delta E(M) = \frac{1}{4} \int_{M} (\Delta H + 2H(H^2-K)+2 e^{-4
\alpha} (A\bar{Z}_3 ^2 + \bar{A}
 Z_3^2)) \varphi d \mu.
 \end{equation}
\end{theorem}

{\it Proof.} By \cite{BT}, we have to compute
\begin{equation}
\label{euler}
\delta E = \frac{1}{4} \left(\delta \int_{M} H^2 d
\mu - \delta \int_{M} \frac{1}{4} n_3 ^{2}\right) d \mu.
\end{equation}
Given parameters $x$ and $y$ on $M$ we mean by $r_1$ and $r_2$ the
derivatives $\frac{\partial r}{\partial x}$ and $\frac{\partial
r}{\partial y}$.

The Gauss--Weingarten equations reads
$$
\nabla_{j} r_i = \Gamma
_{ij} ^k r_k  + h_{ij} n, \ \ \
 \nabla_i n = - h_i^j r_j  = -g^{kj}h_{ki} r_j,
$$
where $g_{ij}$ and $h_{ij}$ are the first and second fundamental
forms of $M$. By the definition of the mean curvature, we have
$$
\delta d\mu = -2H\varphi d\mu.
$$

a) Let us compute
$$
\delta \int H^2d\mu = 2 \int H \delta H d\mu + \int H^2 \delta
d\mu = 2 \int H (\delta H) d\mu - 2 \int H^3 \varphi d\mu.
$$
We have
\begin{equation}
\label{deltaH} 2\delta H = \delta 2H =  \delta (g^{ij} h_{ij}) =
(\delta g^{ij}) h_{ij} + g^{ij} \delta{h_{ij}}
\end{equation}
and the Gauss--Weingarten equations imply that $\delta
g_{ij}=\delta \langle r_i, r_j \rangle = - 2 \varphi h_{ij}$.
Since $0 = \delta \left( g_{ij} g^{jk}\right) =g_{ij}  \delta
g^{jk} + g^{jk} \delta
 g_{ij} = g_{ij}  \delta g^{jk} - 2\varphi h_{ij} g^{jk}$, we have
$$
\delta g^{ij} = 2\varphi g^{jk} h_k^i.
$$

Let us compute $\delta h_{ij}$ which equals $\delta h_{ij} =
\delta \langle \nabla_j r_i, n \rangle =  \langle \nabla_j r_i,
\delta n \rangle + \langle \delta \nabla_j r_i, n \rangle$. By the
Gauss--Weingarten equations, we have $\delta n = - g^{ij}\varphi_j
r_i$ which implies $\langle \nabla_j r_i, \delta n> = -
\Gamma_{ij} ^k \varphi_{k}$. We also have $\delta \nabla_j r_i =
\nabla_{\partial t} \nabla_j r_i = \nabla_j \nabla_{\partial t}
r_i + (\nabla_{\partial t} \nabla_j r_i - \nabla_j
\nabla_{\partial t} r_i) = \nabla_{j} \nabla_{i} (\varphi) n +
\varphi R(r_j,n)r_i$ from which by straightforward computations we
derive that $\langle \nabla_j \nabla_i \varphi n, n \rangle =
\varphi_{ij} - \varphi h_i^k h_{kj}$. Combining the previous
computations we obtain
$$
\delta h_{ij} = -\Gamma _{ij} ^k \varphi _k + \varphi _{ij} -
\varphi h_i ^k h_{kj} + \varphi \langle R(r_j,n)r_i,n \rangle.
$$

Substituting the derived formulas for $\delta g^{ij}$ and $\delta
h_{ij}$ into \eqref{deltaH} we conclude that
$$
2 \delta H = g^{ij} (\varphi_{ij} - \Gamma_{ij}^k \varphi_k) +
\varphi g^{jk}h_k^i h_{ij} + \varphi g^{ij}
\langle R(r_j,n)r_i,n\rangle =
$$
$$
= \Delta \varphi + \varphi h_k ^i h_i ^k + \varphi g^{ij} \langle
R(r_j,n)r_i,n \rangle
$$
where $\Delta$ is the Laplace--Beltrami operator on the surface.
Since $h_k^i h_i^k = \mathrm{Tr}\, h^2 = k_1^2 + k_2^2 = (k_1 +
k_2)^2 -2k_1 k_2 = 4H^2 - 2K$, we rewrite the previous formula as
$$
2 \delta H = \Delta \varphi + (4H^2 - 2K) \varphi + g^{ij} \langle
R(r_j,n)r_i,n \rangle \varphi.
$$

Therefore, by using the equality $\int_M (\Delta \varphi) H  d
\mu= \int _M (\Delta  H) \varphi d \mu $ we derive
$$
\delta \int_{M} H^2 d \mu  =
\int_M (\Delta H + 2H(H^2 - K) \varphi + H g^{ij} \langle
R(r_j,n)r_i,n \rangle) \varphi d\mu.
$$

Assuming that the coordinates $x,y$ are curvilinear orthogonal:
$g_{12}=0$, we have
$$
g^{ij} \langle R(r_j,n)r_i,n \rangle = g^{11} \langle
R(r_1,n)r_1,n \rangle + g^{22} \langle R(r_2,n)r_2,n \rangle =
$$
$$
=
\widehat{K}(r_1,n) + \widehat{K}(r_2,n)
$$
where $\widehat{K}(u,v)$ is the sectional curvature of the ambient
space along the plane spanned by $u$ and $v$.

Let us specialize the formula for $\delta \int H^2 d\mu$ for the
case of surfaces in $\nil$. In this case the sectional curvature
depends only on $n_3$ and equals $\frac{1}{4}-n_3^2$ (see, for
instance, \cite{BT}) and we have $\widehat{K}(r_1,n) +
\widehat{K}(r_2,n) = {n_3}^2 - \frac{1}{2}$ which implies
\begin{equation}
\label{deltaH2}
\delta \int_M H^2 d \mu = \int_M (\Delta H + 2H(H^2 -K) + H(n_3 ^2
- \frac{1}{2}))\varphi d \mu.
\end{equation}

b) Let us compute
$$
\delta \int_{M} {n_3}^2 d \mu = \int 2n_3
\delta n_3 d \mu + \int_{M} {n_3}^2 \delta d \mu.
$$
Therewith we assume that $z=x+iy$ is the conformal parameter on
the surface and the metric takes the form $e^{2\alpha}dzd\bar{z}$.

Since $\langle n, \delta e_3 \rangle = \langle n, \nabla_{\varphi
n} e_3 \rangle = \langle n, \varphi (\frac{1}{2} n_2 e_1 -
\frac{1}{2}n_1 e_2) \rangle = \frac{1}{2} \varphi (n_2 n_1 -n_1
n_2 ) = 0$, \footnote{Here we use the formulas for the Levi-Civita
connection on $\nil$ exposed, for instance, in \cite{BT}.} we have
$\delta n_3 = \delta \langle n,e_3 \rangle = \langle \delta n, e_3
\rangle$ which is equal to $\langle \delta n, e_3 \rangle =
\langle -g^{ij} \varphi_j r_i, e_3 \rangle = -2e^{-2\alpha}
\langle \varphi_z r_{\overline z} + \varphi_{\overline z} r_z ,
e_3 \rangle$.

Thus we compute
$$
2 \int_M n_3 \delta n_3 d \mu = -4 \int_M n_3 \langle \varphi_z
r_{\overline z} + \varphi_{\overline z} r_z,e_3 \rangle dx \wedge
dy =
$$
$$
= 4 \int_{M} ((n_3 \langle r_{\overline{z}},e_3 \rangle)_z + (n_3
\langle r_z, e_3 \rangle)_{\overline{z}})\varphi dx \wedge dy.
$$
By the Weierstrass representation formulas (see \S \ref{sec2}), we
have $\langle r^{-1} r_z, e_3 \rangle =Z_3$ and, since the metric
is left-invariant, we conclude that
$$
\int_M n_3 \delta n_3 d \mu = 2 \int_M (n_3 (\partial
\overline{Z}_3  + \overline{\partial} Z_3) + (\langle
\nabla_{\partial z} n, e_3 \rangle + \langle n, \nabla_{\partial
z} e_3 \rangle)
 \overline {Z}_3 +
$$
$$
+  (\langle \nabla_{\partial \overline z} n, e_3 \rangle + \langle
n, \nabla_{\partial \overline z} e_3 \rangle)
 Z_3 ) \varphi dx \wedge dy
$$
Since $\partial \overline{Z}_3 + \overline{\partial } Z_3 = H n_3
e^{2\alpha} $ (see \cite{BT}), it follows from \eqref{conformal}
that $ \langle \nabla_{\partial z} n, e_3 \rangle = -H Z_3 - 2 A
e^{-2\alpha} \overline{Z}_3$. From the formulas for the
Levi-Civita connection on $\nil$ it also follows that $ \langle n,
\nabla_{\partial z} e_3 \rangle = \langle n , \nabla_{Z_1 e_1 +
Z_2 e_2 + Z_3 e_3}  e_3 \rangle =  \langle n_1 e_1 + n_2 e_2 + n_3
e_3, \frac{1}{2} Z_2 e_1 - \frac{1}{2} Z_1 e_2 \rangle =
\frac{1}{2} (n_1 Z_2 - n_2 Z_1) = \frac{i}{2} Z_3$. Substituting
these formulas into the formula for $\int_M n_3 \delta n_3 d \mu$
we obtain
$$
\int_M n_3 \langle \delta n, e_3\rangle  d \mu = 2 \int_M H n_3 ^2
\varphi d \mu  +
$$
$$
2 \int_M (- 2 H |Z_3|^2 -2 A \bar{Z}_3 ^2 e^{-2\alpha} - 2
  \bar{A}Z_3 ^2 e^{-2\alpha}) \varphi dx \wedge dy
$$
Since $4 |Z_3|^2 = e^{2 \alpha} (1 -n_3 ^2)$ (see
\eqref{metrics}), we have
$$
\int_M n_3 \langle \delta n, e_3 \rangle d \mu = \frac{1}{2}
\int_M (6 H n_3 ^2 -2H - 8 e^{-4 \alpha} (A \bar{Z}_3 ^2 + \bar{A}
Z_3 ^2)) \varphi d \mu,
$$
and finally derive
\begin{equation}
\label{deltan2} \delta \int_M n_3^2 d \mu = \int_M (4Hn_3 ^2 -2H -
8 e^{-4\alpha} (A \bar{Z}_3^2 + \bar{A} Z_3 ^2)) \varphi d\mu.
\end{equation}

Now by substituting \eqref{deltaH2} and \eqref{deltan2} into \eqref{euler} we
prove the theorem.

\subsection{The isoperimetric problem in $S^2 \times \R$ and a
certain Will\-more-type functional} \label{sec6.2}

\begin{proposition}
\label{product}
For a nonminimal constant mean curvature sphere
$M$ in $S^2 \times \R$ we have
$$
\int_M (H^2 + \widehat{K}+1)d\mu = 16\pi.
$$
\end{proposition}

{\it Proof}. By \cite{AR} we know that each cmc sphere is a sphere
of revolution for $H \neq 0$ and a spherical section for $H=0$. By
\cite{Pedrosa} every constant mean curvature sphere of revolution
$S_H$, with $H \neq 0$, is generated by the curve curve $\gamma_H
= S_H / O(2)$ meeting the following equations $B = S^2 \times \R /
O(2) = \{ (x,y) : x \in \R, y \in [0,\pi] \}$:
$$
\frac{dx}{ds} = \cos \sigma, \ \ \frac{dy}{ds} = \sin \sigma, \ \
\frac{d \sigma }{ds} =h +  \cot y \cos \sigma,
$$
where $\sigma$
is the angle between $\gamma$ and the $x$-axis. Thus $\widehat K =
\sin^2 \sigma $ and we have
\begin{equation}
\label{prod1} \int_{S_H} \widehat{K} d\mu  = 4 \pi \left(2 -
\frac{h^2}{\sqrt{1+h^2}} \ln \frac {\sqrt{1+h^2} +1}{\sqrt{1+h^2}
-1}\right).
\end{equation}
By \cite{Pedrosa}, the area of $S_H$ equals
\begin{equation}
\label{prod2} A(S_H) = \int_M d \mu = 4\pi \left(\frac{2}{1+h^2} +
\frac{h^2}{(1+h^2)^{3/2}} \ln \frac {\sqrt{1+h^2} +1}{\sqrt{1+h^2}
-1}\right).
\end{equation}
Combining \eqref{prod1} and \eqref{prod2} we obtain the proof of
the proposition.

As it was proved by Pedrosa \cite{Pedrosa} the isoperimetric
problem for $S^2 \times \R$ is solved by domains bounded by cmc
spheres for small volumes $d$ and by products cylinders $S^2
\times \left[0,\frac{d}{4\pi}\right]$ for large volumes $d$ with
one point $d_0$ of transition from one topological class of
solutions to another.

Proposition \ref{product} one more time demonstrates that a
certain generalization of the Willmore functional of the type
$(H^2 + \alpha \widehat{K}+\beta)d\mu$ respects the isoperimetric
surfaces of spherical topology by attaining on them some constant
value which is even probably the minimum of the functional.

\end{document}